\documentclass[11pt]{article}
\usepackage{amscd,amsmath, amssymb, fancyhdr}

\newcommand{\version}{version 2.0,\ \  Apr. 29, 2014}

\renewcommand{\leftmark}{{\scriptsize LCK metrics obtained from pseudoconvex shells}}

\newcommand{\la}{\lambda}
\newcommand{\al}{\alpha}
\newcommand{\be}{\beta}
\newcommand{\ga}{\gamma}

\newcommand{\f}{\varphi}

\newcommand{\Ll}{\operatorname{Lie}}

\newcommand{\CC}{\mathbb{C}}

\newcommand{\RR}{\mathbb{R}}

\numberwithin{equation}{section}

\def\eqref#1{(\ref{#1})}

\newcommand{\arrow}{{\:\longrightarrow\:}}
\newcommand{\Z}{{\mathbb Z}}
\newcommand{\C}{{\mathbb C}}
\newcommand{\R}{{\mathbb R}}

\newcommand{\N}{{\mathbb N}}

\def\1{\sqrt{-1}\:}
\newcommand{\restrict}[1]{{\left|_{{\phantom{|}\!\!}_{#1}}\right.}}
\newcommand{\cntrct}                % contraction with a vector field
{\hspace{2pt}\raisebox{1pt}{\text{$\lrcorner$}}\hspace{2pt}}

\newcommand{\calo}{{\cal O}}
\newcommand{\cac}{{\cal C}}

\renewcommand{\bar}{\overline}
\renewcommand{\phi}{\varphi}
\renewcommand{\epsilon}{\varepsilon}

\renewcommand{\leq}{\leqslant}

\newcommand{\const}{\operatorname{\text{\sf const}}}

\newcommand{\Aut}{\operatorname{Aut}}

\renewcommand{\phi}{\varphi}

\newcounter{Mycounter}[section]
\newcounter{lemma}[section]
\renewcommand{\thelemma}{{Lemma \thesection.\arabic{lemma}}}
\newcommand{\lemma}{%
     \setcounter{lemma}{\value{Mycounter}}
     \refstepcounter{lemma}
     \stepcounter{Mycounter}
     {\noindent \bf \thelemma:\ }}

\newcounter{claim}[section]
\newcounter{sublemma}[section]
\newcounter{corollary}[section]
\newcounter{theorem}[section]
\renewcommand{\thetheorem}{{Theorem \thesection.\arabic{theorem}}}
\newcommand{\theorem}{%
     \setcounter{theorem}{\value{Mycounter}}
     \refstepcounter{theorem}
     \stepcounter{Mycounter}
     {\noindent \bf \thetheorem:\ }}

\newcounter{conjecture}[section]
\newcounter{proposition}[section]
\renewcommand{\theproposition}
       {{Proposition \thesection.\arabic{proposition}}}
\newcommand{\proposition}{%
     \setcounter{proposition}{\value{Mycounter}}
     \refstepcounter{proposition}
     \stepcounter{Mycounter}
     {\noindent \bf \theproposition:\ }}

\newcounter{definition}[section]
\renewcommand{\thedefinition}
       {{Definition~\thesection.\arabic{definition}}}
\newcommand{\definition}{%
     \setcounter{definition}{\value{Mycounter}}
     \refstepcounter{definition}
     \stepcounter{Mycounter}
     {\noindent \bf \thedefinition:\ }}

\newcounter{example}[section]
\renewcommand{\theexample}{{Example \thesection.\arabic{example}}}
\newcommand{\example}{%
     \setcounter{example}{\value{Mycounter}}
     \refstepcounter{example}
     \stepcounter{Mycounter}
     {\noindent \bf \theexample:\ }}

\newcounter{remark}[section]
\renewcommand{\theremark}{{Remark \thesection.\arabic{remark}}}
\newcommand{\remark}{%
     \setcounter{remark}{\value{Mycounter}}
     \refstepcounter{remark}
     \stepcounter{Mycounter}
     {\noindent \bf \theremark:\ }}

\newcounter{problem}[section]
\newcounter{question}[section]
\makeatletter

\@addtoreset{equation}{section}
\@addtoreset{footnote}{section}
\makeatother

\def\blacksquare{\hbox{\vrule width 5pt height 5pt depth 0pt}}
\def\endproof{\blacksquare}

\begin{document}
%%%%%%%%%%%%%%%%%%%%%%%%%%%%%%%%%%%%%%%%%%%%%%%%%%%%%%%%%%%%
\begin{center}
{\LARGE\bf Locally conformally K\"ahler metrics obtained from pseudoconvex shells
}\\
%%%%%%%%%%%%%%%%%%%%%%%%%%%%%%%%%%%%%%%%%%%%%%%%%%%%%%%%%%%%

\hfill

Liviu Ornea\footnote{Partially supported by CNCS UEFISCDI, project
number PN-II-ID-PCE-2011-3-0118.}
and Misha Verbitsky\footnote{Partially supported by the
 RFBR grant 10-01-93113-NCNIL-a, Science Foundation of 
the SU-HSE award No. 10-09-0015 and AG Laboratory HSE, 
RF government grant, ag. 11.G34.31.0023\\


}

\end{center}

\hfill

%%%%%%%%%%%%%%%%%%%%%%%%%%%%%%%%%%%%%%%%%%%%%%%%
{\small
\hspace{0.15\linewidth}
\begin{minipage}[t]{0.7\linewidth}
{\bf Abstract} \\
A locally conformally K\"ahler (LCK) manifold is a 
complex manifold $M$ admitting a K\"ahler covering $\tilde
M$, such that its monodromy acts on this covering by homotheties.
A compact LCK manifold is called {\bf LCK with potential} 
if its covering admits an authomorphic K\"ahler potential.
It is known that in this case $\tilde M$ is an
algebraic cone, that is, the set of all non-zero vectors in
the total space of an anti-ample line bundle over a
projective orbifold. We start with an algebraic
cone $C$, and show that the set of K\"ahler metrics
with potential which could arise from an LCK 
structure is in bijective correspondence with the set
of pseudoconvex shells, that is, pseudoconvex hypersurfaces
in $C$ meeting each orbit of the associated
$\R^{>0}$-action exactly once. This is used to 
produce explicit LCK and Vaisman metrics on 
Hopf manifolds, generalizing earlier work
by Gauduchon-Ornea and Kamishima-Ornea.
\end{minipage}
}
%%%%%%%%%%%%%%%%%%%%%%%%%%%%%%%%%%%%%%%%%%%%%%%%

\tableofcontents

%%%%%%%%%%%%%%%%%%%%%%%%%%%%%%%%%%%%%%%%%%%%%%%%%%%%%%%%%

\section{Introduction}

%%%%%%%%%%%%%%%%%%%%%%%%%%%%%%%%%%%%%%%%%%%%%%%%%%%%%%%%%

%%%%%%%%%%%%%%%%%%%%%%%%%%%%%%%%%%%%%%%%%%%%%%%%%%%%%%%%%
\subsection{Constructions of LCK metrics}
%%%%%%%%%%%%%%%%%%%%%%%%%%%%%%%%%%%%%%%%%%%%%%%%%%%%%%%%%

\definition
A {\bf locally conformally K\"ahler} (LCK) manifold is a 
complex manifold $M$,  $\dim_\C M >1$, admitting
a K\"ahler covering $(\tilde M, \tilde \omega )$, with
the deck transform group acting on $(\tilde M, \tilde \omega)$
by holomorphic homotheties.

\hfill

For equivalent definitions and examples, see
\cite{_Dragomir_Ornea_} or \cite{_OV:survey_}.

\hfill

A (linear) Hopf manifold $H$ is a quotient of $\C^n\backslash 0$
by an action of $\Z$ generated by a linear map
$A:\; \C^n \arrow \C^n$, with all eigenvalues 
satisfying $|\alpha_i|<1$. It is easy to see that 
the Hopf manifold is diffeomorphic to $S^1\times S^{2n-1}$.
Since $b_1(H)=1$ (an odd number), $H$ is non-K\"ahler.
In fact, this manifold is probably the earliest example
of a non-K\"ahler complex manifold known in mathematics.

However, any Hopf manifold is locally conformally
K\"ahler. This observation originated in works of Izu
Vaisman of 1970-ies; Vaisman defined and studied a strictly smaller class
of manifolds, called by him ``generalized Hopf''. Now these
manifolds are known as ``Vaisman manifolds'', because the name ``generalized Hopf manifold'' was already used by Brieskorn and van de Ven (see \cite{bvv}) for some products of homotopy spheres which do not bear Vaisman's structure. Besides, now it is known that not all Hopf manifolds belong to this smaller class.

These constructions are  non-elementary.
In fact, even the existence of locally conformally
K\"ahler metrics on many Hopf surfaces is quite 
non-trivial. A clean but complicated argument
was given in \cite{go}, where such metrics were  constructed
in dimension 2; when the map $A$ is diagonal, the constructed metric is Vaisman. For non-diagonal
$A$ the construction was much less explicit; in fact,
the LCK metric on these Hopf surfaces was obtained
only by deformation.

In the present paper we use the formalism
of ``locally conformally K\"ahler metrics with potential'',
that we built in previous papers, {\em e.g.} \cite{_OV:LCK_pot_}, to obtain an explicit, 
computation-free and extremely simple
construction of LCK metrics
on manifolds which are obtained as $\Z$-quotients
of algebraic varieties.
This gives, among other things, the first 
explicit ({\em i.e.} not by deformations)  construction of an LCK metric on a
non-diagonal Hopf manifolds.

%%%%%%%%%%%%%%%%%%%%%%%%%%%%%%%%%%%%%%%%%%%%%%%%%%%%%%%%%
\subsection{LCK manifolds}
%%%%%%%%%%%%%%%%%%%%%%%%%%%%%%%%%%%%%%%%%%%%%%%%%%%%%%%%%

This section contains the definitions to be used in
the paper. Unless otherwise stated, we only refer to
compact, connected manifolds (although the definitions
work also for noncompact manifolds).

\hfill

\definition
A complex Hermitian manifold
 $(M,J,g)$ is {\bf locally
conformally K\"ahler} if its fundamental
two-form $\omega:=g\circ J$ satisfies
$$d\omega=\theta\wedge\omega, \quad d\theta=0;$$
here $\theta$ is called {\bf the Lee form.}
This definition is equivalent to the one given above.

A particular subclass of LCK manifolds is described in the following:

\hfill

\definition
A {\bf Vaisman manifold} is a LCK manifold whose Lee form is parallel with respect to the Levi Civita connection of $g$.

\hfill

Compact Vaisman manifolds are equipped with
a Riemannian submersion (a suspension in fact) to the circle with 
the fibers isometric to a Sasakian manifold $N$, see
\cite{ov03} (and see \cite{bg} for an introduction to Sasakian geometry). Their universal coverings are Riemannian
cones $N\times \R^{>0}$ on which the deck transform group,
isomorphic to $\Z$, acts by $(x,t)\mapsto (\varphi(x),
qt)$, where $\varphi$ is a Sasakian automorphism of $N$
and $q\in\N$. The diagonal Hopf manifold is a typical
example, see \S \ref{surv}. On the other hand, it is
known, \cite{belgun}, that non-diagonal Hopf surfaces can
never be Vaisman (although they are LCK, see \cite{go,
  _OV:LCK_pot_}).

A still wider subclass is the following:

\hfill 

\definition
An LCK manifold $M$ which admits a K\"ahler covering 
$(\tilde M, \tilde \omega)$ with the
K\"ahler form $\tilde \omega$ having a global, automorphic potential is
called {\bf LCK manifold with potential}.
Here, by {\bf an automorphic potential} we understand
a function $\psi:\; \tilde M \arrow \R$ satisfying
$dd^c \psi=\tilde \omega$, with the monodromy of $\tilde
M$ mapping $\psi$ to $\const \cdot \psi$.

\hfill 

All Vaisman manifolds are LCK with potential (given by the squared norm of the Lee form with respect to the K\"ahler metric). As for Vaisman manifolds, the monodromy of LCK with potential manifolds is isomorphic to $\Z$. And hence this subclass is strict, as shown by the example of the LCK Inoue surfaces and of the LCK Oeljeklaus-Toma manifolds (see \cite{ot}).

\hfill
 
In this paper we shall be concerned with the {\bf linear Hopf
manifolds.} These are quotients of $\C^n\setminus\{0\}$ by
the cyclic group generated by a linear operator with
eigenvalues strictly smaller than $1$ in absolute value. It is known that
Hopf manifolds are LCK with potential, see
\cite{_OV:LCK_pot_}, and Vaisman if the operator is
diagonal.

%%%%%%%%%%%%%%%%%%%%%%%%%%%%%%%%%%%%%%%%%%%%%%%%%%%%%%%%%
\subsection{Survey of literature}\label{surv}
%%%%%%%%%%%%%%%%%%%%%%%%%%%%%%%%%%%%%%%%%%%%%%%%%%%%%%%%%

There are several papers where explicit constructions
of LCK metric on diagonal Hopf manifolds appear. 

The first one is \cite{vai1}, where the metric (therein named after W. Boothby) $\displaystyle\frac{\sum dz_i\otimes dz_i}{|\sum z_i\bar z_i|^2}$ was considered on $\C^n\setminus\{0\}/\langle z_i\mapsto 2z_i\rangle$. 

More than twenty years took to pass from operators $A=\al\cdot I_n$, $\al\in\CC$, to diagonal operators with complex non-equal eigenvalues. In \cite{go}, a LCK metric was constructed on diagonal
Hopf surfaces $H_{\al,\be}:=\C^2\setminus \{0\}/\langle(u,v)\mapsto(\al u, \be v)\rangle$. The construction is based on
finding a K\"ahler potential on $\C^2\setminus \{0\}$ in terms
of $\al, \be$, but the formula is only implicit. This
procedure was generalized in \cite{belgun}. The
construction of Vaisman metrics in the present paper
can also be considered as a generalization of \cite{go}
to arbitrary dimensions.

In \cite{ko} a construction was done for LCK metrics on 
$\C^n\setminus\{0\}/\langle z_i\mapsto \al_iz_i\rangle$, starting from a deformation of the standard Sasakian structure of $S^{2n-1}$ according to the $S^1$ action with weights $\al_i$ (cf. also \cite[Section 3]{go}). The paper \cite{ko} also contains a very useful criterion to decide when a conformal class of LCK metrics on a complex manifold contains a Vaisman representative, in terms of the existence of a holomorphic complex flow which lifts to a non-trivial flow of homotheties of the K\"ahler covering.    

A different construction on the same manifold, writing explicitly a K\"ahler
potential on $\C^n\setminus \{0\}$,   
appeared in \cite{_Verbitsky_vanishing_} and since then it
was cited in almost all our subsequent papers. However, as
observed by Matei Toma and Ryushi Goto, that metric is 
singular.

To correct this error we provide here a general
construction of LCK metrics on Hopf manifolds. Our
approach works for LCK manifolds with potential, giving
a complete list of LCK metrics with potential in terms
of a pseudoconvex shells in the covering (\ref{_pseudoco_shell_Definition_}).

For the time being, we don't know whether is
it possible or not to write a formula for a potential
for an LCK metric on a Hopf manifold; in the present
paper, as well as in \cite{go}, the potential is written as
a solution of a certain differential equation.

%%%%%%%%%%%%%%%%%%%%%%%%%%%%%%%%%%%%%%%%%%%%%%%%%%%%%%%%%

\section{Algebraic cones and LCK manifolds with potential}

%%%%%%%%%%%%%%%%%%%%%%%%%%%%%%%%%%%%%%%%%%%%%%%%%%%%%%%%%

%%%%%%%%%%%%%%%%%%%%%%%%%%%%%%%%%%%%%%%%%%%%%%%%%%%%%%%%%
\subsection{Algebraic cones}
%%%%%%%%%%%%%%%%%%%%%%%%%%%%%%%%%%%%%%%%%%%%%%%%%%%%%%%%%

%%%%%%%%%%%%%%%%%%%%%%%%%%%%%%%%%%%%%%%%%%%%%%%%%%%%%%
\definition\label{_algebra_cone_intro_Definition_}
{\bf A closed algebraic cone}  is an
affine variety $C$ admitting a $\C^*$-action $\tau$
with a unique fixed point $x_0$ (called {\bf the origin}), which satisfies
the following.

\begin{description}
\item[1.] $C$ is smooth outside of $x_0$.
\item[2.] $\tau$ acts on the Zariski tangent
space $T_{x_0}\cac$ diagonally, with all eigenvalues
$|\alpha_i|<1$.
\end{description}

\noindent {\bf An open algebraic cone} is a closed algebraic cone with the origin removed: $C \setminus \{x_0\}$.

\hfill

%%%%%%%%%%%%%%%%%%%%%%%%%%%%%%%%%%%%%%%%%%%%%%%%%%%%%%%%
\definition
Let $X$ be a projective
orbifold, and let $L$ be an ample line bundle on $X$. 
Assume that the total space of $L$ is smooth
outside of the zero divisor. The {\bf algebraic
cone $\cac(X,L)$ of $X,L$} is the total space
of non-zero vectors in $L^*$.
A {\bf cone structure} on $\cac(X,L)$ is the $\C^*$-action
arising this way (by fiberwise multiplication).

\hfill

In \cite[Section 4]{_OV:Sasakian_on_CR_}, it was shown that
any open algebraic cone $C$ is isomorphic to $\cac(X,L)$,
for appropriate $X$ and $L$. This was shown
by the following argument. Given the algebraic cone $C$, one obtains $X$ as the quotient
of $C$ by $\C^*$, and then the
cone $C$ is naturally identified with the total
space of a principal $\C^*$-bundle $L_1$. Ampleness
of this bundle follows, because the corresponding
closed algebraic cone $C_c$ is an affine variety,
and algebraic functions on $C_c$ are identified with the
section of the line bundle $L$ associated with $L_1$.

\hfill

%%%%%%%%%%%%%%%%%%%%%%%%%%%%%%%%%%%%%%%%%%%%%%%%%%%
\definition
Let $\gamma$ be an automorphism of a closed algebraic cone.
It is called {\bf a holomorphic contraction} if for any 
compact subset $K\subset C$, and any open neighbourhood
$U$ of the origin, there exists a number
$N$ sufficiently big such that $\gamma^N(K)\subset U$.

\hfill

%%%%%%%%%%%%%%%%%%%%%%%%%%%%%%%%%%%%%%%%%%%%%%%%%%%
\definition
Let $C$ be a closed algebraic cone, and $\rho:\; \R^{>0}\arrow \Aut(C)$
a $\R^{>0}$-action. We say that $\R^{>0}$ acts {\bf by holomorphic
contractions}, if $\rho(t)$ is a holomorphic contraction
for all $t<1$.

\hfill

%%%%%%%%%%%%%%%%%%%%%%%%%%%%%%%%%%%%%%%%%%%%%%%%%%%
\example
Let $C=\C^n=\cac(\C P^{n-1}, \calo(1))$.
Then any linear automorphism of $C$ with all
eigenvalues $|\alpha_i|<1$ acts on $C$
by holomorphic contractions.

\hfill

%%%%%%%%%%%%%%%%%%%%%%%%%%%%%%%%%%%%%%%%%%%%%%%%%%%
\example
Let $\rho:=\tau\restrict{\R^{>0}}$ 
be the action of $\R^{>0}$ on an algebraic cone
provided by the cone structure, $\R^{>0}\subset \C^*$.
Since $\rho$ acts on the tangent space
$T_cC$ to the origin with eigenvalues
smaller than 1, it acts on $C$ by holomorphic contractions
(\cite[Theorem 3.3]{_OV:LCK_pot_}).

\hfill

%%%%%%%%%%%%%%%%%%%%%%%%%%%%%%%%%%%%%%%%%%%%%%%%%%%
\remark 
As shown in \cite[Theorem 3.3]{_OV:LCK_pot_},
the quotient of an algebraic cone by a contraction
is an LCK manifold with potential, and, conversely,
any  LCK manifold with potential is obtained 
by taking the quotient of an open algebraic cone by 
a holomorphic contraction. Such a contraction, being 
{\em a priori} a $\Z$--action, can be extended 
to a $\R^{>0}$--action by holomorphic contractions.

\hfill

Further on, we use the following version of this result.

\hfill

%%%%%%%%%%%%%%%%%%%%%%%%%%%%%%%%%%%%%%%%%%%%%%%%%%%
\theorem
Let $M$ be a locally conformally K\"ahler manifold with potential.
Then $\tilde M$, as a complex manifold,
 is isomorphic to an open algebraic cone $C$, equipped 
with an action $\rho$ of $\R^{>0}$ by holomorphic
contractions, and the quotient $\tilde M/\langle \rho(2^n)\rangle$
is isomorphic (as a complex manifold) to $M$.

\hfill

{\bf Proof:} In \cite[Theorem 2.1]{_OV:top_}, it is shown that
$M$ can be deformed into a Vaisman manifold. From its
proof it is apparent that this deformation preserves
$\tilde M$ (in fact, only the $\Z$-action is deformed). 
Therefore, $\tilde M$ is a covering of a Vaisman manifold.
Then, it is an algebraic cone, as follows from
\cite[Proposition 4.6]{_OV:LCK_immer_}.
\endproof

%%%%%%%%%%%%%%%%%%%%%%%%%%%%%%%%%%%%%%%%%%%%%%%%%%%%%%%%%
\subsection{CR-geometry and Sasakian manifolds}
%%%%%%%%%%%%%%%%%%%%%%%%%%%%%%%%%%%%%%%%%%%%%%%%%%%%%%%%%

In this subsection,
we introduce the Sasakian manifolds and some 
related notions of CR-geometry. We follow 
\cite{_OV:Sasakian_on_CR_}.

\hfill

%%%%%%%%%%%%%%%%%%%%%%%%%%%%%%%%%%%%%%%%%%%%%%%%%%%%%%%%%
\definition
 A {\bf CR-structure} (Cauchy-Riemann structure) on a manifold 
$M$ is a subbundle $H\subset T M \otimes \C$ of the complexified tangent bundle, which is closed
under  commutator:
$[H, H ] \subset  H$ 
and satisfies $H \cap \bar H = 0$.

A function $f:M\rightarrow \C$ is {\bf CR-holomorphic} if $D_Vf=0$ for any vector field $V\in \bar H$.

\hfill

On a CR manifold $(M,H)$, the bundle $H\oplus \bar H$ is preserved by complex conjugation and hence it is obtained as a complexification of a real subbundle $H_\R$. Then $I_H:=-\sqrt{-1}\text{Id}_{\bar H}$ defines a complex structure on $H_\R$ and $H$ is its $\sqrt{-1}$-eigenspace of its extension to the complexification $H_\R \otimes \C$.

If $\text{codim}_{TM}H_\R=1$,  and the {\em Frobenius tensor} $L:H_\R\times H_\R\rightarrow TM/H_\R$, $L(X,Y)=[X,Y] \mod H_\R$ is nondegenerate, then $(M,H)$ is a {\bf CR contact manifold} and $H_\R$ is its contact structure (or distribution). In this context, $L$ is called the {\bf Levi form}.

As $L$ vanishes on $H$ and $\bar H$, $L$ is $(1,1)$ with respect to $I_H$. 

\hfill

%%%%%%%%%%%%%%%%%%%%%%%%%%%%%%%%%%%%%%%%%%%%%%%%%%%%%%%%%
\definition
A contact CR-manifold $(M, H_\R , I_ H )$ is called {\bf
  pseudo-convex} 
if the Levi form is positive or negative, depending on the choice of orientation. If
this form is also sign-definite, then $(M, H_\R , I_H )$ is called {\bf strictly pseudoconvex}.

\hfill

%%%%%%%%%%%%%%%%%%%%%%%%%%%%%%%%%%%%%%%%%%%%%%%%%%%%%%%%%
\definition
Let $S$ be a CR-manifold.
A CR-holomorphic vector field $v\in TS$ 
is called {\bf transversal} if it is transversal to 
the CR-distribution $H_\R\subset TS$.

\hfill

%%%%%%%%%%%%%%%%%%%%%%%%%%%%%%%%%%%%%%%%%%%%%%%%%%%%%%%%%%%%
\theorem\label{_S^1_equiv_Sasakian_Theorem_} \cite[Theorem 1.2]{_OV:Sasakian_on_CR_} 
Let $M$ be a compact pseudoconvex contact CR-manifold.
Then the following conditions are equivalent.
\begin{description}
\item[(i)] $M$ admits a Sasakian metric, compatible
with the CR-structure.
\item[(ii)] $M$ admits a proper, transversal
CR-holomorphic $S^1$-action.
\item[(iii)] $M$ admits a nowhere degenerate, transversal
CR-holomorphic vector field.
\end{description}

%{\bf Proof:} \cite[Theorem 1.2]{_OV:Sasakian_on_CR_}. \endproof

\hfill

%%%%%%%%%%%%%%%%%%%%%%%%%%%%%%%%%%%%%%%%%%%%%%%%%%%%%%%%%%%%
\theorem\label{_Sasa_unique_Theorem_} \cite[Theorem 1.3]{_OV:Sasakian_on_CR_} 
Let $M$ be a compact, strictly pseudoconvex CR-manifold
admitting a proper, transversal CR--ho\-lo\-mor\-phic $S^1$--ac\-tion.
Then $M$ admits a unique (up to an automorphism)
$S^1$-invariant CR-embedding into an algebraic cone $(\cac,\tau)$.
Moreover, a Sasakian metric on $M$ can be induced
from a K\"ahler metric $\tilde\omega$ on this 
cone, which is automorphic in the following sense:
for some constant $c>1$, one has $\tau(t)^*\tilde\omega=|t|^c\tilde\omega$.

\hfill

%{\bf Proof:} \cite[Theorem 1.3]{_OV:Sasakian_on_CR_}. \endproof

%%%%%%%%%%%%%%%%%%%%%%%%%%%%%%%%%%%%%%%%%%%%%%%%%%%%%%%%%
\subsection{Pseudoconvex shells in algebraic cones}
%%%%%%%%%%%%%%%%%%%%%%%%%%%%%%%%%%%%%%%%%%%%%%%%%%%%%%%%%

\hfill

%%%%%%%%%%%%%%%%%%%%%%%%%%%%%%%%%%%%%%%%%%%%%%%%%%%%%%%%%
\definition\label{_pseudoco_shell_Definition_}
Let $C$ be an algebraic cone, equipped 
with an action $\rho$ of $\R^{>0}$ by holomorphic
contractions. A {\bf pseudoconvex shell}
in $C$ is a strictly pseudoconvex submanifold
in $C$, intersecting each orbit
of $\rho$ exactly once.

\hfill

%%%%%%%%%%%%%%%%%%%%%%%%%%%%%%%%%%%%%%%%%%%%%%%%%%%%%%%%%
\remark
Please note that the action of $\rho$ may bear no
relation to the cone action $\tau: \C^* \arrow \Aut(C)$.

\hfill

%%%%%%%%%%%%%%%%%%%%%%%%%%%%%%%%%%%%%%%%%%%%%%%%%%%%%%%%%
\theorem\label{_shell_bijective_to_potentials_Theorem_}
Let $M=C/\langle\rho(q)\rangle$ where $(C,\tau)$ be an  algebraic cone,
equipped with the action $\rho$ of $\R^{>0}$ by holomorphic
contractions and $q>1$. Let $\vec{r}$ be the infinitesimal generator of $\rho$ and let $S$ be a pseudoconvex shell in $C$. 
Then for each $\la\in\R$ there exists a unique function $\phi_\la$ such that $\Ll_{\vec{r}}\phi_\la=\la\phi_\la$ and $\phi_\la\restrict{S}=1$. Moreover, such $\f_\la$ is plurisubharmonic for sufficiently big  $\la>\!\!>0$.
Conversely, any LCK manifold with potential admits a metric obtained this way.

\hfill

{\bf Proof:} For each $\rho$-orbit and each $\rho$-equivariant potential $\phi$, one has:
$$\rho(t)\cdot \phi_\la=e^{t\la} \phi_\la, \qquad t\in \R^{>0}.$$
Let $S$ be a pseudoconvex shell in $C$. Then $S\times \R^{>0}\stackrel{\sim}{\longrightarrow} C$, as orbits intersect the shell only once. Hence for any $s\in S$ and for any $t\in\R^{>0}$ we have:
\begin{equation}\label{rot}
\phi_\la(\rho(t)\cdot s)=e^{\la t},
\end{equation}
and this equation uniquely defines $\phi_\la$. The problem is to prove that $dd^c\f_\la>0$ and this is not automatic for $\f_\la$, but it holds for some power of it.

Now let $B:=e^{\R\vec{r}} \cdot (TS\cap I(TS)) \subset TC$ be the subbundle obtained by 
translating $TS\cap I(TS)$ with all $e^{t\vec{r}}$. Then, by construction, $dd^c\f_\la \restrict{B}$ is the Levi form of $B$ and hence it is positive definite.

It will now suffice to show that $dd^c\f_{2a\la}\restrict{S}=dd^c\f_\la ^{2a}\restrict{S}>0$ for sufficiently big $a$. But 
$$dd^c\f_\la^{2a}=\f_\la^{2a-2}\big(2a\f_\la \cdot dd^c\f_\la +2a(2a-1)d\f_\la\wedge d^c\f_\la \big).$$
As the shell $S$ is compact, the result is implied by the following elementary linear algebra lemma:

\hfill

\lemma\label{linear} Let $h_1$, $h_2$ be pseudo-Hermitian forms on a complex vector space $V$ and let $W\subset V$ be a codimension 1 subspace. Assume that $h_1\restrict{W}$ and $h_2\restrict{V/W}$ are strictly positive 
(that is, positive definite), and $h_2\restrict{W}=0$. 
Then there exists $u_0\in\R$, depending continuously on $h_1, h_2$, such that $h_u:=h_1+uh_2$ is positive definite for any $u>u_0$.

\hfill

The direct part of the Theorem now follows by applying \ref{linear} (whose proof we postpone) to $V=TM$, $W=B$, $h_1=\f_\la  dd^c\f_\la $, $h_2=d\f_\la\wedge d^c\f_\la $.

For the converse, let $\f_\la$ be any automorphic potential, thus satisfying $(\gamma^k)^*=e^\la\f_\la$, and let $\vec{r}$ be the holomorphic vector field which is the logarithm of the monodromy action. Let then $\rho(t)=e^{-t\la}\big(e^t\vec{r}\big)^*$ the corresponding endomorphism of $\cac ^{\infty}(M)$. Then $\rho(k+t)(\f_\la)=\rho(t)\f_\la$ and hence the orbit of $\rho$ through $\f$ is compact. We then average $\rho(t)\f$ on $\RR$ and obtain a $\rho(t)$-invariant K\"ahler potential $\f_{\la 0}$. This $\f_{\la 0}$ is obtained from a pseudoconvex shell $S=\f_{\la 0}^{-1}(1)$ and from $\vec{r}$ as in the direct part of the Theorem.

\hfill

Let us now give the proof of \ref{linear}. For simplicity, we work in the real setting, and we consider $h_1, h_2$ as bilinear symmetric forms. Let $y\in V$ be a vector such that $h_2(y,y)=1$. Then any $x\in V$ can be written as $x=ay+z$, for some $z\in W$. This translates to:
\begin{equation*}
h_u(x,x)=ua^2+a^2h_1(y,y)+h_1(z,z)+2ah_1(z,y)
\end{equation*}
which we view as a polynomial in $a$. This one is positive definite for all $a$ if and only if
\begin{equation}\label{**}
(h_1(z,y))^2-(u+h_1(y,y))\cdot h_1(z,z)<0.
\end{equation}
Choose $y'\in W$ such that $h_1(z,y')=h_1(z,y)$ for all $z\in W$ and let $u>u_0:=h_1(y',y')-h_1(y,y)$. Then \eqref{**} becomes
$$(h_1(z,y'))^2-h_1(y',y')h_1(z,z)<0,$$
which is satisfied by Cauchy-Buniakovski-Schwarz inequality, as $h_1$ is positive definite on $W$.

\hfill \endproof

\example Let $A$ be a linear operator on $\C^n$ with
eigenvalues of absolute values strictly smaller than
$1$. Let $C=\C^n\setminus \{0\}$ and let $\rho(t)=e^{t\log
A}$. Take a sphere $S=S^{2n-1}\subset \C^n$; it is easy to
see that $S$ is a pseudoconvex shell. Applying
\ref{_shell_bijective_to_potentials_Theorem_}, we obtain 
an automorphic potential $\phi$ on $\C^n\backslash 0$,
giving an LCK metric on $M=(\C^n\backslash 0)/\langle
A\rangle$.

 \hfill
 
\remark
When $n=2$ and $A$ is diagonal,  the same potential
was obtained in \cite{go}.

In particular, we recover the result proven in 
\cite{go,ko, _OV:LCK_pot_} that {\em all Hopf manifolds
  $(\C^n\setminus \{0\})/\langle A\rangle$ are LCK}.

%%%%%%%%%%%%%%%%%%%%%%%%%%%%%%%%%%%%%%%%%%%%%%%%%%%%%%%%%
\subsection{Vaisman metrics and pseudoconvex shells}

\remark Vaisman manifolds are LCK with potential and hence
they have canonical pseudoconvex shells (levels of the
potential).

\hfill

%%%%%%%%%%%%%%%%%%%%%%%%%%%%%%%%%%%%%%%%%%%%%%%%%%%%%%%%%
\definition
Let $(M,I, g)$ be a Vaisman manifold, let 
$(C, \rho)$ be the associated algebraic cone,
equipped with the action $\rho$ of $\R^{>0}$ by holomorphic
contractions, and let $S$ be its pseudoconvex shell.
The {\bf Reeb field} of $M$ is the  CR-holomorphic vector
field $I\theta^\sharp\in TS$ 
obtained from $\rho$ by complex conjugation.

\hfill

%%%%%%%%%%%%%%%%%%%%%%%%%%%%%%%%%%%%%%%%%%%%%%%%%%%%%%%%%
\remark
For Vaisman manifolds, the Reeb field is always transversal
(\cite{_OV:Sasakian_on_CR_}).

\hfill

\proposition
Let $v\in TS$ be a CR-holomorphic vector field on $S$,
where $S$ is a pseudoconvex shell in an algebraic cone $C$. 
Then $v$ can be uniquely extended to a holomorphic
vector field on the whole of $C$. \footnote{This is called
{\bf the holomorphic extension} of  $v$.}

\hfill

{\bf Proof:} Let   $\calo_S$, respectively
$\calo_{{S^\circ}}$ be the ring of CR-holomorphic
functions on $S$, respectively on the interior ${S^\circ}$
of the shell $S$. By the solution of the Neumann problem,
$L^2(\calo_S)=L^2(\calo_{{S^\circ}})$ and hence, if we
restrict to bounded functions, $\calo_{{S^\circ}}=\calo_S$
as rings. As vector fields on $S$ and ${S^\circ}$ are
derivations of the above rings, the result follows.
\endproof

\hfill

%%%%%%%%%%%%%%%%%%%%%%%%%%%%%%%%%%%%%%%%%%%%%%%%%%%%%%%%%
\theorem\label{explicit_construction}
Let $(M,I,g)$ be an LCK manifold obtained (as in
\ref{_shell_bijective_to_potentials_Theorem_})
from an algebraic cone $C$ and a pseudoconvex shell $S$, 
and let $\gamma:\; \Z \arrow \Aut(C)$ be the deck transform map.
Then the Hermitian
manifold $(M,I,g)$ is conformally equivalent to
a Vaisman one if and only if 
$S$ admits a transversal CR-holomorphic
vector field $\xi$, such that its holomorphic extension 
to $C$ is $\gamma(1)$--invariant and $\exp(-I\xi)\cdot\ga(1)$ preserves $S$.

\hfill

{\bf Proof:}  Suppose $(M,I,g)$ is Vaisman. Then
$\xi=I\theta^\sharp$ is an isometry of the LCK
metric. The K\"ahler metric on $C$ is $dd^c\f$, with $\f$
given by the equation
$$e^{t\theta^\sharp}(S)=\f^{-1}(t),$$
where $S$ is the level set of $\f$.
Then $S$ is Sasakian and $I\theta^\sharp$ is the Reeb
field of the underlying contact structure, hence it is
transversal by definition. By construction,
$M=C/\langle\gamma(1)\rangle$, and the extension of $\xi$ to
$C$ is the  lift  to $C$ and is $\gamma(1)$ invariant by
definition. We have checked all the conditions of
\ref{explicit_construction}, except the last one:
$\exp(-I\xi)\cdot\ga(1)(S)=S$.

The Lee field $-I\xi=\theta^\sharp$  acts by
homotheties on the potential $\f$, hence
$\Ll_{\exp(-I\xi)}\f=c\f$ for some contant $c$. The
monodromy map $\ga(1)$ also acts by homotheties on $\f$:
$\Ll_{\rho(1)}\f=c'\f$ for another constant $c'$. What we
have to do is to homothetically modify the initial metric
$g$ such that $c$ becomes $1/c'$; in this case
$\exp(-I\xi)\cdot\ga(1)$ will preserve the potential $\f$
and hence will preserve  $S$.

Conversely, the vector field $\xi$ is tangent to $S$ and
its flow $\exp(t\xi)$ acts by contractions, hence defining
a metric on the cone $C$ over $S$. As $\gamma(1)$
maps a shell $S_\la$ to $S_{\const \cdot\lambda}$, 
if follows that $\gamma(1)$ acts
by homotheties on $C$. This means that the complex flow
generated on $M=C/\langle\gamma(1)\rangle$ by $\xi$ and
$I\xi$ lifts to a  flow of non-trivial homotheties on the
cone. By \cite{ko} this ensures the existence of a Vaisman
metric in the conformal class of $g$.
\endproof

%%%%%%%%%%%%%%%%%%%%%%%%%%%%%%%%%%%%%%%%%%%%%%%%%%%%%%%%%
\subsection{Examples and {\em erratum}}
%%%%%%%%%%%%%%%%%%%%%%%%%%%%%%%%%%%%%%%%%%%%%%%%%%%%%%%%%

As an application, we add an {\em erratum} to several papers 
where we have given a wrong formula for an LCK-metric
on diagonal Hopf manifolds ({\em e.g.} \cite{_Verbitsky_vanishing_}, \cite{ornea}), and give a general
and almost explicit construction of a Vaisman metric
on any diagonal Hopf manifold. This constructions originates in the ones in \cite{go, belgun, ko}, but unifies them and presents them in a 
much more synthetic and transparent way.

\hfill

We apply \ref{explicit_construction}. The data to start
with are the cone $C$, the shell $S$, the vector field
$\xi$ and the monodromy $\ga$. 

Let $A=\mathrm{diag}(\al_1,\cdots,\al_n)$, with
$0<|\al_1|\leq |\al_2|\leq\cdots\leq|\al_n|<1$. Denote by
$A_{|\cdot|}$ the matrix
$\mathrm{diag}(|\al_1|,\cdots,|\al_n|)$
(in the same basis).

Let $C$ be $\C^n\setminus \{0\}$. As a shell $S$, we take the
sphere $S^{2n-1}$, but one can take for $S$ the boundary of any
strictly pseudoconvex body containing 0 and 
satisfying $A_{|\cdot|}A^{-1}(S)=S$.

The transversal CR-holomorphic vector field $\xi$ is 
\begin{equation}\label{def_xi}
\xi=\sqrt{-1}\mathrm{Re} \log A=\sqrt{-1}\log\mathrm{diag}(|\al_1|,\ldots,|\al_n|)
\end{equation}
and the monodromy map $\ga$ is given by $\ga(z)=A\cdot z$, $z\in C$. 

To apply  \ref{explicit_construction} we need to verify that $\xi$ is $\gamma(1)$--invariant and $\exp(-I\xi)\cdot\ga(1)$ preserves $S$. An easy computation shows that $\exp(-I\xi)\cdot\ga(1)$ acts as $A_{|\cdot|}A^{-1}=\mathrm{diag}(\frac{\al_1}{|\al_1|},\cdots,\frac{\al_n}{|\al_n|})$ and thus it preserves the norm of vectors, hence preserving the spheres. Finally, as the action of $\ga(1)$ is linear, the $\gamma(1)$--invariance of $\xi$ amounts to $A\cdot\xi_z=\xi_{A\cdot z}$ which is immediate from \eqref{def_xi}.

This gives an explicit construction of Vaisman structure on  diagonal Hopf manifolds. On the other hand, \ref{_shell_bijective_to_potentials_Theorem_} provides a rather explicit construction of LCK metrics on non-diagonal Hopf manifolds, avoiding the argument by deformations in \cite{_OV:LCK_pot_}.

\hfill

\noindent{\bf Acknowledgments:} We thank Matei Toma and
Ryushi Goto for drawing our attention to the wrong formula
which originated this paper, and to Max Planck Institute for Mathematics   
(Bonn), where this paper was finished, for the excellent research environment.

\hfill

{\small

\noindent {\sc Liviu Ornea\\
University of Bucharest, Faculty of Mathematics, \\14
Academiei str., 70109 Bucharest, Romania. \emph{and}\\
Institute of Mathematics ``Simion Stoilow" of the Romanian Academy,\\
21, Calea Grivitei Street
010702-Bucharest, Romania }\\
\tt Liviu.Ornea@imar.ro, \ \ lornea@gta.math.unibuc.ro

\hfill

\noindent {\sc Misha Verbitsky\\
{\sc Laboratory of Algebraic Geometry, \\
National Research University HSE,\\
7 Vavilova Str. Moscow, Russia, 117312}\\
\tt  verbit@mccme.ru,\ \  verbit@verbit.ru\\
}}

\end{document}